\documentclass[12pt]{article}
\usepackage{amsmath, amsfonts, amsthm, amssymb, epsfig}

\begin{document}
\title{{\bf A family of diophantine equations of the form
    $x^4+2nx^2y^2+my^4=z^2$ with no solutions in $({\mathbb Z}^+)^3$}}
\author{Konstantine Zelator\\
Department of Mathematics\\
College of Arts and Sciences\\
Mail Stop 942\\
The University of Toledo\\
Toledo, OH  43606-3390\\
e-mails: konstantine-zelator@yahoo.com\\
\ \hspace*{.75in}konstantine.zelator@utoledo.edu}

\maketitle

\section{Introduction}

In this work we present a family of diophantine equations of the form

\begin{equation}
x^4+2nx^2y^2+my^4 = z^2 \label{E1}
\end{equation}

\noindent with no nontrivial solutions.

This is done in Section 3, where the theorem in this paper, Theorem 1, and its
proof are presented.  The approach is elementary and uses only congruence
arguments as well as decent.  It is branched proof, with some of the branches
leading to contradictions via congruence arguments.  Two of the proof's
branches lead to contradictions via a decent argument.  Also in the proof, we
make use of the well-known parametric formulas that describe all the solutions
in $({\mathbb Z}^+)^3$ to the diophantine equation $x^2 + \ell \cdot y^2 =
z^2$, $\ell$ a positive integer.  These formulas are found in Section 2.  In
Section 4, we present a sampling of numerical examples. That is, a listing of
combinations of integers $n$ and $m$ in (\ref{E1}), which satisfy the
hypothesis of the theorem. 

The paper concludes with Section 5, wherein we offer a brief historical
commentary on diophantine equations of the form $ax^4+bx^2y^2+cy^4=dz^2$.
Investigations of these types of diophantine equations span a time interval of
nearly 400 years, not to go back any further in time.  We mention some of the
results found in the literature, including more recent developments (of the
last 70 years) on the subject involving the usage of local methods as well as
the association of such equations with elliptic curves.

\section{An auxiliary diophantine equation: \\
$x^2+\ell \cdot y^2 = z^2$}

For a given positive integer $\ell$, the solution set (subset of $({\mathbb
  Z}^+)^3$) of the diophantine equation $x^2+\ell y^2 = z^2$, can be
  parametrically described by the formulas,

$$
x = \dfrac{d(\rho_1k^2-\rho_2\lambda^2)}{2},\ y = dk\lambda,\ \ z =
\dfrac{d(\rho_1k^2+\rho_2\lambda^2)}{2}
$$

\noindent where the parameters $d,k,\lambda$ are positive integers such that
$(k,\lambda)=1$; and the positive integers $\rho_1,\rho_2$ are divisors of
$\ell$ such that $\rho_1\rho_2 = \ell$.  Obviously, if we require that
$(x,y)=1$, then all the solutions in $({\mathbb Z}^+)^3$ can be parametrically
described as follows:

\begin{equation}
\left\{\begin{array}{l}
x = \dfrac{d(\rho_1k^2-\rho_2 \lambda^2)}{2},\ y = dk\lambda,\ \ z =
\dfrac{d(\rho_1k^2 + \rho_2 \lambda^2)}{2},\\
\\
{\rm with}\ d,k,\lambda,\rho_1,\rho_2 \in {\mathbb Z}^+\ {\rm such\ that}\
(k,\lambda) = 1,\ \rho_1\rho_2 = \ell \\
\\
{\rm and\ with}\ d = 1\ {\rm or}\ 2.\ {\rm Also},\ \rho_1 k^2-\rho_2\lambda^2
> 0. \end{array}\right\}\label{E2}
\end{equation}

These parametric formulas are well known in the literature and can be found in
reference \cite{1}, (pages 420-421).  A derivation of them can also be found
in \cite{2}.

\section{The theorem and its proof}

\noindent {\bf Theorem 1:}  Suppose that $n$ is a positive integer, $p$ an
odd prime, and such that either

$$\begin{array}{lcll}
n \equiv 0\ ({\rm mod}\ 4) & {\rm and} & p \equiv 3\ ({\rm mod}\ 8); & {\rm or\
  alternatively},\\
\\
n \equiv 2\ ({\rm mod}\ 4) & {\rm and} & p \equiv 7\ ({\rm mod}\ 8)
\end{array}
$$

\noindent In addition to the above, assume that one of the following
hypotheses holds:  

\noindent either 

\begin{enumerate}
\item[(i)]  $n^2-p>0$ and the positive integer $m = n^2 - p$ is a prime, or

\item[(ii)]  $n^2 - p < 0$ and the positive integer $N = -m = -(n^2-p)$ is a
  prime.
\end{enumerate}

\noindent Then, the diophantine equation $x^4+2nx^2y^2 +my^4 = z^2$ has no
solution in $({\mathbb Z}^+)^3$.

\vspace{.15in}

\noindent {\bf Proof:}  If equation (\ref{E1}) has a solution, then let
$(X_0,Y_0,Z_0)$ be a solution with the product $X_0Y_0$ being least.  Let
$\delta = (X_0,Y_0)$, so that $X_0 = \delta x_0, \ Y_0 = \delta y_0$, for
positive integers $x_o,y_0,\delta$ such that $(x_0,y_0)=1$.  Then, (\ref{E1})
 implies $\delta^4 \mid Z^2_0$ and so $\delta^2 \mid Z_0$; and by
putting $Z_0 = \delta z_0$ for some $z_0 \in {\mathbb Z}^+$ we obtain 

\begin{equation}
x^4_0 +2nx^2_0 y^2_0+my^4_0 = z^2_0 \label{E3}
\end{equation}

\noindent By (\ref{E3}), the triple $(x_0,y_0,z_0)$ is a solution to equation
(\ref{E1}).  Thus, by the minimality of the product $X_0 Y_0$ it
follows that $\delta = 1,\ X_0=x_0,\ Y_0=y_0,\ Z_0 = z_0$.

Since $x_0$ and $y_0$ are relatively prime, there are three possibilities:
$$
x_0 \equiv y_0 \equiv 1\ ({\rm mod}\ 2);\ x_0 \equiv 0\ {\rm and}\ y_0 \equiv
1 \ ({\rm mod}\ 2); \ {\rm or}\ x_0 \equiv 1\ {\rm and}\ y_0 \equiv 0\ ({\rm
  mod}\ 2).
$$

\noindent If $x_0$ and $y_0$ are both odd, consider equation (\ref{E3})
modulo 4.  Since $x^2_0 \equiv y^2_0 \equiv 1\ ({\rm mod}\ 4)$, in this case,
(\ref{E3}) implies $1+2n+m = z^2_0\ ({\rm mod}\ 4)$.  By the Theorem's
hypothesis, $2n \equiv 0\ ({\rm mod}\ 4)$ and $m=n^2-p$; we obtain $1-p \equiv
z^2_0\ ({\rm mod}\ 4)$, which gives $z^2_0 \equiv 2\ ({\rm mod}\ 4)$ in view
of $p\equiv 3\ ({\rm mod}\ 4)$, an impossibility.  

Next consider the second possibility.  The combination $x_0$ being even and
$y_0$ odd.  By (\ref{E3}), since $m$ is odd, we see that $z_0$ must be odd as
well.  Consider (\ref{E3}) modulo 8.  In view of $y^2_0 \equiv z^2_0 \equiv 1\
({\rm mod} \ 8)$, (\ref{E3}) implies  $m \equiv 1\ ({\rm mod}\ 8)$, a
contradiction since by hypothesis:
$$m=n^2-p \equiv 0 - 3 \equiv 5\ ({\rm mod}\ 8),\  {\rm if}\ n \equiv 0\ ({\rm
  mod}\ 4) \  {\rm and}\ p \equiv 3\ ({\rm mod}\ 8),$$

\noindent  while also,
$$m = n^2-p \equiv 4-7\equiv 5\ ({\rm mod}\ 8),\ {\rm if}\ n \equiv 2\ ({\rm
  mod}\ 4) \ {\rm and}\ p \equiv 7\ ({\rm mod}\ 8).
$$

We conclude that $x_0$ must be odd and $y_0$ even.  Also, it is clear from
(\ref{E3}) that since $(x_0,y_0)=1$, $y_0$ and $z_0$ must be relatively prime
as well; and $z_0$ must be odd.  Therefore,

\begin{equation}
\left\{ \begin{array}{l}x^4_0 + 2nx^2_0y^2_0 + my^4_0 = z^2_0\\
\\
x_0 \equiv z_0 \equiv 1\ ({\rm mod}\ 2),\ y_0 \equiv 0\ ({\rm mod}\ 2)\\
\\
(x_0,y_0) =1=(y_0,z_0) \end{array} \right\} \label{E4}
\end{equation}

Now we use the hypothesis that $m = n^2 - p$.  An algebraic manipulation of
the equation in (\ref{E4}) leads to,

\begin{equation}
\begin{array}{rcl}
(x^2_0 +ny^2_0)^2 - z^2_0 & = & py^4_0;\\
\\
\left[ (x^2_0 + ny^2_0) + z_0\right] \left[(x^2_0 +ny^2_0) - z_0\right]& =&
py^4_0 \end{array}\label{E5} \end{equation}

According to the conditions in (\ref{E4}) both $(x^2_0 + ny^2_0)$ and $z_0$
are odd integers,  but they are also coprime.  Indeed, if a prime $q \neq p$
were a common divisor of theirs, then by (\ref{E5}) it would also divide $y_0$
and therefore,  $x_0$ as well, violating $(x_0,y_0) =1$.  If $p$ divided both
$(x^2_0 +ny^2_0)$ and $z_0$, then $p^2$ would divide the left-hand side of
(\ref{E5}), and thus $p$ would divide $y_0$.  Hence, it would divide $x_0$,
contrary once more to $(x_0,y_0)=1$.  We conclude that

\begin{equation}
(x^2_0 + ny^2_0, z_0) = 1 \label{E6}
\end{equation}

\noindent Moreover, the sum of any two odd integers is congruent to $0\ ({\rm
  mod}\ 4)$ and their difference to $2\ ({\rm mod}\ 4)$; or vice-versa.
  Combining this observation with (\ref{E6}) leads to,

\begin{equation}
\left\{ \begin{array}{l}
x^2_0 +ny^2_0 + z_0 = 2 \delta_1\\
\\
x^2_0 + ny^2_0 - z_0 = 2\delta_2\\
\\
{\rm for}\ \delta_1,\delta_2\in {\mathbb Z}^+,\ {\rm with}\
(\delta_1,\delta_2) = 1\ {\rm and}\ \delta_1+\delta_2 \equiv 1\ ({\rm mod}\ 2)
\end{array}\right\}.\label{E7} \end{equation}

Adding the two equations in (\ref{E7}) yields,

\begin{equation} x^2_0 +ny^2_0 = \delta_1 + \delta_2. \label{E8} 
\end{equation}

According to (\ref{E7}), $\delta_1$ must be even and $\delta_2$ odd; or
vice-versa.  Given that the rest of the proof rests on (\ref{E8}) and that
(\ref{E8}) is symmetric in $\delta_1$ and $\delta_2$.  There is no need to
distinguish between two cases, they lead to the same contradictions.
Accordingly, assume that $\delta_1$ is even and $\delta_2$ is odd.

If we combine (\ref{E7}) with (\ref{E5}), we see that since $p$ is an odd
prime, there are precisely two possibilities expressed in (9) below:

\vspace{.15in}

 \fbox{\parbox{4.5in}{$\begin{array}{lllll}
{\rm Either} & 2 \delta_1 = 8py^4_1 & {\rm and} & 2\delta_2=2y^4_2 & \hspace{1.0in}{\rm (9a)}\\
{\rm or} & 2 \delta_1 = 8y^4_1 & {\rm and} & 2\delta_2 = 2py^4_2
&\hspace{1.0in}  {\rm (9b)}
\end{array}$
\newline for positive integers $y_1,y_2,$ such that $(y_1,y_2) =1$  \newline
and $y_2 \equiv 1\
({\rm mod}\ 2)$. \newline
Note that in either case, we have from (\ref{E5}), $2y_1y_2 = y_0$.}} \hfill
(9)

\vspace{.15in}

\setcounter{equation}{9}

\noindent {\bf Case 1:}  Assume possibility (9b) in (9) to hold.  Then by
combining (9b) with (\ref{E8}) gives

$$
x^2_0 + ny^2_0 = 4y^4_1 + py^2_2,
$$

\noindent which is impossible modulo 4, since by (\ref{E4}) we have $x^2_0
+ny^2_0 \equiv 1\ ({\rm mod}\ 4)$, while $4y^4_1 +py^4_2 \equiv p \equiv 3\
({\rm mod}\ 4)$, in view of the hypothesis of the theorem.

\vspace{.15in}

\noindent {\bf Case 2:}  Assume possibility (9a) to be the case in (9).

\vspace{.15in}

\noindent {\it Subcase 2(i):}  Assume hypothesis (i) in the theorem, which means that
the integer $n^2-p=m$ is positive and a prime.  By combining (9a) with
(\ref{E8}) and $2y_1y_2 = y_0$ in (9) we obtain

\begin{equation}
x^2_0 + (n^2-p)\cdot (2y^2_1)^2 = (y^2_2 -2ny^2_1)^2 \label{E10}
\end{equation}

According to (\ref{E10}), the triple $(x_0,2y^2_1,\left|y^2_2 -
2ny^2_1\right|)$ is a positive integer solution to the diophantine equation
$x^2+\ell y^2=z^2$, with $\ell = n^2-p$.  Also note that $(x_0,2y^2_1)=1$, by
virtue of the fact that $(x_0,y_0)= 1$ in (\ref{E4}) and $y_0 = 2y_1y_2$ in
(9).  Therefore, by (\ref{E2}), we must have

\vspace{.15in}

\noindent \parbox{5.0in}{$\left\{ \begin{array}{l}2y^2_1 = dk\lambda,\
 \left|y^2_2-2ny^2_1\right| = \dfrac{d(\rho_1k^2+\rho_2\lambda^2)}{2},\ {\rm for\ positive}\\
{\rm integers}\  d,k,\lambda,\rho_1,\rho_2;\ {\rm such \ that}\ (k,\lambda)=1,\ \rho_1\rho_2 = n^2-p\\
 {\rm and\ with}\ d = 1\ {\rm or}\ 2\end{array}\right\}$ }\hfill (10a)

\vspace{.15in}

The possibility $d=1$ is easily ruled out by the fact that $\rho_1$ and
$\rho_2$ are both odd (since $m=n^2-p$ is odd); and the fact that
$(k,\lambda)=1$.  Indeed, the first equation (10a) implies, if $d=1$, that $k$
and $\lambda$ must have different parities.  But, then the integer $\rho_1k^2
+ \rho_2\lambda^2$ would be odd, instead of even as the second equation in
(10a) requires.  Thus, $d=2$ which yields, by (10a)

\vspace{.15in}

\parbox{5.0in}{$\begin{array}{l} y^2_1=k\lambda,\ \left|y^2_2 -2ny^2_1\right|
    = \rho_1k^2+\rho_2\lambda^2;\ {\rm or\ equivalently}\\
\\
\left\{ y^2_1 =k\lambda,\ y^2_2-2ny^2_1 = \pm
    (\rho_1k^2+\rho_2\lambda^2)\right\}\end{array}$} \hfill (10b)

\vspace{.15in}

The first equation in (10b) implies, since $(k,\lambda)=1$, that $k=k^2_1\ {\rm and}\ \lambda = \lambda^2_1;\ {\rm for\ some}\ \lambda_1,k_1
\in {\mathbb Z}^+,\ {\rm with}\ (k_1,\lambda_1)=1.$

Accordingly (10b) gives,

\vspace{.15in}

\hspace*{1.25in}$y^2_2 - 2nk^2_1 \lambda^2_1 = \pm (\rho_1k^4_1 + \rho_2
\lambda^4)$ \hfill (10c)

\vspace{.15in}

\noindent since $y_1 = k_1\lambda_1$.

\vspace{.15in}

If the plus sign holds in (10c), we obtain

\vspace{.15in}

\hspace*{1.25in} $y^2_2 = \rho_1k^4_1 + 2nk^2_1\lambda^2_1+\rho_2\lambda^4_1$
\hfill (10d)

\vspace{.15in}

By (10a) we know that $\rho_1\rho_2 = m =n^2-p$.  But $m$ is a prime and so
either $\rho_1=m$ and $\rho_2 = 1$ or vice-versa.  In either case, (10d) shows
that the triple $(k_1,\lambda_1,y_2)$ is a positive integer solution to the
diophantine equation (\ref{E1}).  Compare this solution with the solution
$(x_0,y_0,z_0)$ (see (3)).  We have 

$$
x_0y_0 \geq y_0 = 2y_1y_2 > y_1 = k_1 \lambda_1.
$$

\noindent In short, by (9) $x_0y_0 > k_1\lambda_1$, contradicting the fact that
$x_0y_0$ is least.

If the minus sign holds in (10c),

\vspace{.15in}

\hspace*{1.25in} $y^2_2 = - \rho_1k^4_1 + 2nk^2_1\lambda^2_1 - \rho_2
\lambda^4_1$ \hfill (10e)

\vspace{.15in}

\noindent Again, we use the fact that either $\rho_1 = n^2-p$ and $\rho_2 =1$
or vice-versa.

In either case, $\rho_1+\rho_2 = n^2 - p+1$.  Consider (10e) modulo 4.  If
both $k_1$ and $\lambda_1$ are odd, then $k^2_1 \equiv \lambda^2_1 \equiv 1\
({\rm mod}\ 4)$ and so (10e) implies,
$$
\begin{array}{l}
y^2_2 \equiv - \rho_1 +2n - \rho_2 \equiv -(\rho_1+\rho_2) + 2n \equiv -
(n^2-p+1)+2n\ ({\rm mod}\ 4);\\
\\
y^2_2 \equiv - n^2 + p - 1 + 2n \equiv 2 \ ({\rm mod}\ 4),\end{array}
$$

\noindent since by hypothesis  $p\equiv 3\ ({\rm mod}\ 4)$ and $n$ is
even. Thus, a contradiciton.

If $k_1+\lambda_1 \equiv 1\ ({\rm mod}\ 2)$, again consider (10e) modulo 4.
Given that $\rho_1 = n^2-p$ and $\rho_2 =1$ or vice-versa, and that $k_1$ is
odd and $\lambda_1$ even, or vice-versa.  The four combinations, because of
the symmetry of (10e) reduce to two congruence possibilities:  $y^2_2 \equiv
-1$ or $y^2_2 \equiv -(n^2-p)\ ({\rm mod}\ 4)$,  but $n^2-p \equiv 1\ ({\rm
  mod}\ 4)$, by hypothesis .  Therefore we see that in both cases we arrive at
$y^2_2 \equiv 3\ ({\rm mod}\ 4)$ which is impossible.  This concludes the
proof in subcase (2i).

\vspace{.15in}

\noindent {\it Subcase 2(ii):}  Assume hypothesis (ii) of the theorem.  Then
$n^2-p<0$ and $N=p-n^2$ is a prime.  Combining (\ref{E8}) with (9a) and
$2y_1y_2=y_0$ in (9) leads to

\begin{equation}
x^2_0 = (y^2_2 - 2ny^2_1)^2 + (p-n^2)(2y^2_1)^2 \label{E11}
\end{equation}

By (9) we know that $(y_1,y_2) = 1$ and $y_2$ is odd; which implies that
$(y^2_2-2ny^2_1, 2y^2_1) =1$.  By (\ref{E11}), the triple $
\left(\left|y^2_2-2ny^2_1\right|, 2y^2_1,x_0\right)$ is a positive integer
solution to the diophantine equation $x^2+\ell y^2 = z^2$, with $\ell = p-n^2$;
and with the integers $\left|y^2_2 - 2ny^2_1\right|$ and $2y^2_1$ being relative prime.  Accordingly, by (\ref{E2})  we must have

$$\left| y^2_2 - 2ny^2_1\right| = \dfrac{d(\rho_1k^2-\rho_2\lambda^2)}{2},\
2y^2_1 = dk\lambda;
$$

\begin{equation}
\left\{ \begin{array}{l}
y^2_2 - 2ny^2_1 = \pm \dfrac{d(\rho_1k^2 - \rho_2 \lambda^2)}{2},\ 2y^2_1 =
dk\lambda,\\ 
{\rm for \ positive \ integers} \ d,k,\lambda,\rho_1,\rho_2 \ {\rm such\ 
  that}\\
(k,\lambda)=1,\ \rho_1\rho_2 = p-n^2,\ {\rm and\ with}\ d=1\ {\rm or}\
2\end{array}\right\} \label{E12} \end{equation}

Since we consider (below) all the combinations $\rho_1,\rho_2$ such that
\linebreak $\rho_1\rho_2 = p-n^2$, it follows that the plus or minus possibilities in the first equation of (\ref{E12}) are really the same.  Thus, we may write

\vspace{.15in}

\hspace*{1.0in} $y^2_2 - 2ny^2_1 = \dfrac{d(\rho_1k^2 - \rho_2\lambda^2)}{2},\
2y^2_1 = dk\lambda$ \hfill (12a)

\vspace{.15in}

As we saw in the proof of subcase (ii), the possibility $d=1$ is easily ruled
out.  Indeed, if $d=1$, the first equation in (12a) implies that the integer
$(\rho_1k^2 -\rho_2 \lambda^2)$ must be even.  On the other hand, the second
equation in (12a) implies, since $(k,\lambda)=1$ that $k$ must be odd and
$\lambda$ even; or vice-versa.  But then, by virtue of the fact that
$\rho_1,\rho_2$ are both odd, it follows that $\rho_1k^2-\rho_2 \lambda^2
\equiv 1\ ({\rm mod}\ 2)$, a contradiction.  Thus, $d=2$ in (12a).  We have,

\vspace{.15in}

\hspace*{1.0in} $y^2_2-2ny^2_1 = \rho_1k^2 - \rho_2\lambda^2,\ y^2_1 =
k\lambda$ \hfill (12b)

\vspace{.15in}

Obviously, the second equation in (12b) implies, since $(k,\lambda)=1$, that
$k=k^2_1$ and $\lambda^2_1 = \lambda$ for some $k_1,\lambda_1 \in {\mathbb
  Z}^+$, with $k_1,\lambda_1)=1$.  Using $y_1 = k_1\lambda_1$ as well, we
see that (12b) implies

\vspace{.15in}

\hspace*{1.5in} $y^2_2 = \rho_1k^4_1 + 2nk^2_1\lambda^2_1 - \rho_2\lambda^4_1$
\hfill (12c)

\vspace{.15in}

Since $\rho_1\rho_2 = p-n^2 =$ prime, there are precisely two possibilities.
Either $\rho_1=1,\ \rho_2=p-n^2$ or, alternatively, $\rho_1=p-n^2$ and $\rho_2
=1$.  In the first case, $\rho_1 =1$ and $-\rho_2 = n^2-p=m$; so that by
(12c), $y^2_2 = k^4_1 +2nk^2_1\lambda^2_1 + m\lambda^4_1$, which shows that
the triple $(k_1,\lambda_1,y_2)$ is a positive integer solution to the
initial equation (\ref{E1}).  Compare this solution with the solution
$(x_0,y_0,z_0)$.  We have, $x_0y_0 \geq y_0 = 2y_1y_2 > y_1 = k_1 \lambda_1$,
violating the minimality of the product $x_0y_0$.  Next, assume the next
possibility to take hold, namely $\rho_1=p-n^2$ and $\rho_2 =1$.  Then
equation (12c) implies,

\vspace{.15in}

\hspace*{1.5in} $y^2_2 = (p-n^2)k^4_1 + 2nk^2_1 \lambda^2_1 - \lambda^4_1$
\hfill (12d)

\vspace{.15in}

\noindent Consider (12d) modulo 4:

\vspace{.15in}

\noindent If $k_1 \equiv \lambda_1 \equiv 1\ ({\rm mod}\ 2)$, then (12d) implies
$y^2_2 \equiv p-n^2 + 2n-1\ ({\rm mod}\ 4) \Rightarrow$ (since $n$ is even and
$p\equiv 3\ ({\rm mod}\ 4))\ y^2_2 \equiv 2\ ({\rm mod}\ 4)$, an
impossibility.  

\noindent If $k_1 \equiv 0$ and $\lambda_1 \equiv 1\ ({\rm mod}\ 2)$, (12d)
implies $y^2_2 \equiv -1 \equiv 3\ ({\rm mod}\ 4)$, again impossible.

\noindent Finally, if $k_1$ is odd and $\lambda_1$ even, (12d) implies $y^2_2
\equiv p-n^2\ ({\rm mod}\ 4);$ \linebreak $ y^2_2 \equiv 3\ ({\rm mod}\ 4)$, again an impossibility.

This concludes the proof of subcase (ii) and with it, the proof of the
theorem. \hfill \rule{2mm}{2mm}

\section{Numerical Examples}

\begin{enumerate}
\item[(i)]  Below, we provide a list of all combinations of positive integers
  $n,p,m$; such that both $p$ and $m$ are primes, $m=n^2-p$, and with either
  $n \equiv 0\ ({\rm mod}\ 4)$ and $p\equiv 3\ ({\rm mod}\ 8)$, or alternatively,
  $n \equiv 2\ ({\rm mod}\ 4)$ and $p \equiv 7\ ({\rm mod} \ 8)$.  Under the
  constraint $n \leq 16$, there are 24 such combinations.

$$
\begin{array}{|r|c|c|c|}
\hline
&n & p & m\\ \hline
1) & 4 & 3 & 13\\ \hline
2)  & 4 & 11 & 5 \\ \hline
3) & 6 & 7 & 29 \\ \hline
4) & 6 & 23 & 13 \\ \hline
5) & 6 & 31 & 5 \\ \hline
6) & 8 & 3 & 61 \\ \hline
7) & 8 & 11 & 53 \\ \hline
8) & 8 & 59 & 5 \\ \hline
9) & 10 & 47 & 53 \\ \hline
10) & 10 & 71 & 29 \\ \hline
11) & 12 & 43 & 101 \\ \hline
12) & 12 & 83 & 61 \\ \hline
\end{array} \hspace*{.5in}
\begin{array}{|r|c|c|c|}
\hline
& n & p & m \\ \hline
13) & 12 & 107 & 37 \\ \hline
14) & 12 & 131 & 13 \\ \hline
15) & 12 & 139 & 5 \\ \hline
16) & 14 & 23 & 173\\ \hline
17) & 14 & 47 & 149 \\ \hline
18) & 14 & 167 & 29 \\ \hline
19) & 14 & 191 & 5 \\ \hline
20) & 16 & 59 & 197 \\ \hline
21) & 16 & 83 & 173\\ \hline
22) & 16 & 107 & 149 \\ \hline
23) & 16 & 227 & 29 \\ \hline
24) & 16 & 251 & 5 \\ \hline 
\end{array}
$$

\item[(ii)]  Below, we provide a listing of all combinations of integers $n,
  p, m, N$; such that $n,p, N > 0,\ m<0,\ p$ and $N$ are both primes,
  $N=p-n^2,\ m=-N,$ and with either $n \equiv 0\ ({\rm mod}\ 4)$ and $p \equiv
  3\ ({\rm mod}\ 8)$, or alternatively, with $n \equiv 2 \ ({\rm mod}\ 4)$ and
  $p \equiv 7\ ({\rm mod}\ 8)$. Under the constraint $p \leq 251$, there are
  29 such combinations.

$$\begin{array}{|r|c|c|c|c|}
\hline
& p & n & N & m \\ \hline
1) & 7 & 2 & 3 & -3 \\ \hline
2) & 23 & 2 & 19 & -19 \\ \hline
3) & 47 & 2 & 43 & -43 \\ \hline
4) & 47 & 6 & 11 & -11 \\ \hline
5) & 59 & 4 & 43 & -43 \\ \hline
6) & 67 & 8 & 3 & -3 \\ \hline
7) & 71 & 2 & 67 & -67 \\ \hline
8) & 79 & 2 & 73 & -73 \\ \hline
9) & 79 & 6 & 43 & -43 \\ \hline
10) & 83 & 4 & 67 & - 67 \\ \hline
11) & 83 & 8 & 19 & -19 \\ \hline
12) & 103 & 6 & 67 & -67 \\ \hline
13) & 103 & 10 & 3 & -3 \\ \hline
14) & 107 & 8 & 43 & -43 \\ \hline
15) & 131 & 8 & 67 & - 67 \\ \hline
\end{array} \hspace*{.5in}
\begin{array}{|r|c|c|c|c|} 
\hline
& p & n & N & m \\ \hline
16) & 163 & 12 & 19 & -19 \\ \hline
17) & 167 & 2 & 163 & -163 \\ \hline
18) & 167 & 6 & 131 & -131 \\ \hline
19) & 167 & 10 & 67 & -67 \\ \hline
20) & 179 & 4 & 163 & -163\\ \hline
21) & 199 & 6 & 163 & -163 \\ \hline
22) & 199 & 14 & 3 & -3 \\ \hline
23) & 211 & 12 & 67 & -67 \\ \hline
24) & 227 & 4 & 211 & -211 \\ \hline
25) & 227 & 8 & 163 & -163 \\ \hline
26) & 227 & 12 & 83 & -83 \\ \hline
27) & 239 & 10 & 139 & -139 \\ \hline
28) & 239 & 14 & 43 & -43 \\ \hline
29) & 251 & 12 & 107 & -107 \\ \hline
& & & &  \\ \hline
\end{array}
$$ 

\section{Historical Commentary}

Mathematical research on diophantine equations of the form
\begin{equation}
ax^4+bx^2y^2 + cy^4 = dz^2 \label{E13}
\end{equation}

\noindent dates back to the early 17th century.  The most comprehensive source
of results on such equations in the 300 year-period from the early 17th
century to about 1920, is I. E. Dickson's monumental book {\it History of the
  Theory of Numbers, Vol. II}, (see [1]).

All or almost all results (at least the referenced ones) of that period can be
found in that book.  Various researchers during that time period employed
decent methods to tackle such equations.  Perhaps all the significant results
achieved in that 300-year period can be attributed to about 40-50
investigators.  We list the names of thirty-two of them:

\vspace{.15in}

\noindent Fermat, Frenicle, St. Martin, Genocci, Lagrange, Legendre, Lebesgue,
Euler, Adrain, Gerardin, Aubry, Fauquenbergue, Sucksdorff, Gleizes, Mathieu,
Moret-Blank, Rignaux, Kausler, Fuss, Auric, Realis, Mantel, Desboves, Kramer,
Escott, Thue, Cunningham, Pepin, Lucas, Werebrusov, Carmichael, Pocklington.  

\vspace{.15in}

A detailed account of the results obtained by these mathematicians is given in
[1], pages 615-639.

On the other hand, the last 75 years or so (from the early 1930's to the present)
are marked by the introduction and development of what is known as local
methods as well as the connection/association of equations (\ref{E13}) with
elliptic curves.  In particular, the beginning of the 75 year period (early
thirties) is characterized by a landmark, the Hasse Principle: 

\vspace{.15in}

\noindent If $F \in {\mathbb Z}[x_1,\ldots , x_n ]$ is a homogenous polynomial
of degree $2$, then $F(x_1,\ldots , x_n)$ has a nontrivial solution in
${\mathbb Z}^n$ if, and only if, 

\begin{enumerate}
\item[(a)] it has a nontrivial solution in ${\mathbb R}^n$ and 
\item[(b)] it has a primitive solution modulo $p^k$ , for all primes $p$ and
  exponents $k \geq 1$. 
\end{enumerate}

Here, a solution $(a_1,\ldots , a_n)$ is understood to be nontrivial if at
least one of the $a_i$'s is not zero.  It is primitive if one of the $a_i$'s
is not divisible by $p$.

In 1951, E. Selmer (see \cite{3}), presented an example of a homogenous
polynomial in three variables, and degree $n=3$ which fails the Hasse Principle.This is the equation $3x^3 +4y^3+5z^3 =0$, whose only solution in ${\mathbb
  Z}^3$ is $(0,0,0)$ (so it has no nontrivial solutions).  But it obviously
has nontrivial solutions in ${\mathbb R}^3$; and it has primitive solutions
modulo each prime power.

In their paper W. Aitken and F. Lemmermeyer, (see \cite{4}), show that
equation (\ref{E13}) has a nontrivial solution in ${\mathbb Z}^3$ if, and
only if, the diophantine system (in four variables $u,v,w,z$)

\begin{equation}
\left\{ \begin{array}{ll}
{\rm with}\ b^2-4ac \neq 0, & au^2 + bv^2 + cw^2 = dz^2 \\
\\
{\rm and} \ d\ {\rm squarefree}, & uw = v^2 \end{array}\right\} \label{E14}
\end{equation} 

\noindent has a nontrivial solution in ${\mathbb Z}^4$.  This also holds when
${\mathbb Z}$ is replaced by any ring containing ${\mathbb Z}$.  In
particular, it holds for ${\mathbb R}$.

Furthermore, (\ref{E14})  has a primitive solution modulo $p^k$ if, and only
if, (\ref{E13})  has a primitive solution modulo $p^k$; and $k \geq 2$.  (If
$p$ is not a divisor of $d$, this can be extended to $k=1$.)

In 1940 and 1942 respectively, C.-E Lind and H. Reichardt, (see \cite{5} and
\cite{6}), found another counterexample to the Hasse Principle:  the 
diophantine equation (\ref{E13}) with $a=1,\ b=0,\ c = -17$, and $d = 2$; that
is the equation $x^4-17y^4=2z^2$.

Aitken and Lemmermeyer generalized the Lind and Reichardt example by taking
$a=1,\ b=0,\ c = -q$, such that $q$ is a prime  with $q \equiv 1\ ({\rm mod}\
16),\ d$ is squarefree, $d$ is a nonzero square but not a fourth power modulo
$q$, and $q$ is a fourth power modulo $p$ for every odd prime $p$ dividing
$d$.  Thus, they obtained a family of diophantine equations (\ref{E13}) (or
equivalently, systems (\ref{E14})) which fail the Hasse Principle.  Their
proofs of the nontrivial insolvability (of each member of that family)  in
${\mathbb Z}^3$ only involves quadratic reciprocity arguments.  The harder
part is to give an elementary proof that the above equations have primitive
solutions modulo all prime powers.

Variants of the Hasse Principle, and the manner in which these principles
fail, can be found in a paper by B. Mazur (see \cite{7}).  Also, there is the
seminal work by J. Silverman (see \cite{8}), which provides a comprehensive
study for the links between equations (\ref{E13}) and elliptic curves.

Alongside these developements of the last 75 years, there have been some
results obtained by elementary means only.  For example, A. Wakulitz (see
\cite{9}) has offered an elementary proof that the diophantine equation
$x^4+9x^2y^2 +27y^4 = z^2$ has no solution in $({\mathbb Z}^+)^3 $.  A
corollary of this (in the paper in \cite{9}), is that the equation
$x^3+y^3=2z^3$ has no solution in ${\mathbb Z}^3$ with $x \neq y$ and $z \neq
0$.

\end{enumerate}

\end{document}